\documentclass[12pt]{amsart}
\usepackage{amsfonts}
\usepackage{ifthen}
\usepackage{amsthm}
\usepackage{amsmath}
\usepackage{graphicx}
\usepackage{amscd,amssymb,amsthm}

\newcounter{minutes}
\divide\time by 60
\newcounter{hours}
\multiply\time by 60 \addtocounter{minutes}{-\time}

\title[V. N. Dubinin, M. Vuorinen/Extremal Decomposition problem]
{An extremal decomposition problem for harmonic measure$^\bigstar$}

\author[]{Vladimir N. Dubinin$\dagger$}
\author[]{Matti Vuorinen$\ddagger$}

\address{Institute of Applied Mathematics, Far-Eastern Branch of Russian Academy of Sciences, Vladivostok, Russia} \email{dubinin@iam.dvo.ru}

\address{Department of Mathematics, University of Turku, Turku 20014,
Finland} \email{vuorinen@utu.fi}

\thanks{$\dagger$The research of this author was supported by Far-Eastern Branch of Russian Academy of Sciences, project 09-III-A-01-007
\thanks{$\ddagger$ Supported by the Academy of Finland, project 2600066611}
\thanks{$^\bigstar$ File:~\jobname .tex,
          printed: \number\year-\number\month-\number\day,
          \thehours.\ifnum\theminutes<10{0}\fi\theminutes} }

\newtheorem{theorem}[equation]{Theorem}

\numberwithin{equation}{section}

\begin{document}


\maketitle

\begin{abstract}
Let $E$ be a continuum in the closed unit disk $|z|\le 1$
of the complex $z$-plane which divides the open disk  $|z| < 1$
into $n\ge 2$ pairwise non-intersecting simply connected domains
$D_k,$ such that each of the domains $D_k$ contains
some point $a_k$ on a prescribed circle $|z| = \rho, 0 <\rho <1\, , \, k=1,\dots,n\,. $
It is shown that for some increasing function $\Psi\,,$ independent of $E$
and the choice of the points $a_k,$ the mean value of the harmonic measures
$$ \Psi^{-1}\left[   \frac{1}{n} \sum_{k=1}^{k}  \Psi(  \omega(a_k,E, D_k))\right] $$
is greater than or equal to the harmonic measure $\omega(\rho, E^* , D^*)\,,$ where
$E^* = \{ z: z^n \in [-1,0] \}$ and $D^* =\{  z: |z|<1,   |{\rm arg} z| < \pi/n\} \,.$
This implies, for instance, a solution to a problem of R.W. Barnard, L. Cole, and
A. Yu. Solynin concerning a lower estimate of the quantity
$\inf_{E}  \max_{k=1,\dots,n} \omega(a_k,E, D_k)\,$ for arbitrary points of the circle  $|z| = \rho \,.$ These authors stated this hypothesis in the particular case when the points
are equally distributed on the circle  $|z| = \rho \,.$
\end{abstract}

\section{Introduction and formulation of the result}

Problems concerning extremal decomposition of domains on the complex sphere
go back to the works of Lavrentiev, Goluzin, Nehari, Jenkins
and have a rich history and are more or less directly associated with
several areas of geometric function theory (see, e.g.,  \cite{j3,kuz,s}).
A majority of extremal decompositions is related to estimates for the
products  whose factors are powers of inner radii of pairwise nonoverlapping
domains. For the case of harmonic measure the problem concerning
extremal decomposition was perhaps first formulated in the paper of
R.W. Barnard, L. Cole, and
A. Yu. Solynin \cite{bcs}. Let $E$ be a continuum in the closed unit disk
$\overline{U}, U = \{z: |z| < 1 \}$ and let it divide $U$ into two
subdomains $D_1 \ni \rho$ and $D_2 \ni -\rho, 0 < \rho <1\,.$ Thus
$U \setminus E = D_1 \cup D_2 \,.$ The authors of the paper \cite{bcs} give
the minimum of the sums
$$
\frac{1}{2}( \omega(\rho, E, D_1)  + \omega(-\rho, E, D_2))\,,
$$
taken over the set of all continua $E \subset \overline{U}$ that divide $U$
into a pair of domains as described above.  Here $\omega(z, E, D)$ denotes the harmonic measure of the set $E \cap \partial D$ with respect to the domain $D\,,$ evaluated at the
point $z\,.$ For the properties of the harmonic
measure we refer to \cite[Ch 1.1]{a}.
In the paper \cite{bcs} also a physical motivation of the problem
is given. It is of interest to observe that the extremal configuration in this problem
is not symmetric with respect to the imaginary axis. Among other things the authors of the
paper \cite{bcs} discuss the problem of finding the lower bound
$$
\inf_E  \max_{k=1,\dots,n} \omega(a_k^* , E, D_k)
$$
over all continua $E$ that divide $U$ into $n$ simply connected domains
$$D_k \ni a_k^*= \rho \, {\rm exp} (2 \pi i(k-1)/n), k=1,\dots,n, n\ge 3\,.$$
They conjectured that the extremal configuration of the minmax problem consists
of $n$ circular sectors
$$
D_k^* = \{ z \in U: |{\rm arg} z - 2 \pi(k-1)/n|< \pi/n \}, \,\,\, k=1,\dots,n\,.
$$
and observed that under the additional assumption that the closure of
$\partial D_k \cap U$ is connected for all $k=1,\dots,n,$ this can be proved
applying Jenkins' theory of extremal decomposition \cite{j1,s}(see \cite[p. 246]{bcs}).
A. Yu. Solynin has informed the authors that this result has remained
unpublished.

In this note we prove the following theorem.

\begin{theorem} \label{mainth}
For a given $\rho, 0 < \rho <1,$
let $E$ be a continuum in the closed unit disk $\overline{U}$
dividing the open disk  $U$
into $n\ge 2$ pairwise non-intersecting simply connected domains
$D_k,$ such that each of the domains $D_k$ intersects
the circle $|z| = \rho , k=1,\dots,n\,. $ Then for arbitrary
points $a_k \in D_k$ on the circle $|z| = \rho \,,$ the following
inequality holds
\begin{equation} \label{mainineq}
\frac{1}{n} \sum_{k=1}^n \log \frac{1+ \sin(\pi \omega_k/2)}{1- \sin(\pi \omega_k/2)} \ge -n \log \rho
\end{equation}
where $\omega_k = \omega(a_k,E, D_k), k=1,\dots,n\,.$ Equality holds in {\rm (\ref{mainineq}) } if and only if
 $$E = \{ z: ( e^{i \theta}z)^n \in [-1,0] \},
D_k =\{  z: z e^{i \theta} \in D_k^*\} $$
and $a_k = a_k^* e^{- i \theta}, k=1,\dots,n, $ for some fixed real $\theta\,.$
\end{theorem}

Because the function
$$
\Psi(x) =\log \frac{1+ \sin (\pi \, x/2)}{1- \sin (\pi \, x/2)}
$$
is strictly increasing on the interval $(0, 1),$ the inequality
(\ref{mainineq}) implies a solution to the problem of
R.W. Barnard, L. Cole, and
A. Yu. Solynin \cite{bcs}, and furthermore, the points $a_k$ need
not coincide with $a_k^*, k=1,\dots,n\,.$ Our proof is based on the
result of the first author concerning the product of inner radii
of nonoverlapping domains with respect to "free'' points on the
given circle \cite{d1}. We also make use of an idea of J. Krzyz
about the transition from inner radius to conformal invariants \cite{k}.
If in the case of \cite{k} the Green function were discussed, then
in our case the harmonic measure is in the focus. We note further, that
by the invariance of the harmonic measure under  M\"obius automorphisms
of the disk $U\,,$ the inequality (\ref{mainineq}) holds for all points
$a_k,$ located on an arbitrary hyperbolic circle in $U$ of hyperbolic
radius $ \log (({1+\rho})/({1-\rho})) \,.$




\section{Proof of Theorem \protect{\ref{mainth}}}
Fix an integer $k, 1 \le k \le n\,.$ If all the boundary points of the
domain $D_k$ are lying in the continuum $E,$ then $\omega_k=1$ and the left
side of the inequality (\ref{mainineq}) is $+\infty\,.$ In what follows we
exclude this case for all $k=1,\dots,n\,.$ Let now $z$ be some point of the
boundary $\partial D_k,$ not contained in $E\,.$ Then $z \in T:= \partial U \,.$
In fact, in the opposite case, $z\in U \setminus E\,,$ and therefore $z\in D_j$
for some $j, 1\le j \le n\,.$ Because of the openness of a domain, we see
that $D_j \cap D_k \neq \emptyset,$ which contradicts the hypothesis of the
theorem. Further, because the set $E$ is closed, on the circle $T$ there
is an open arc $\alpha$ containing the point $z$ and not containing any points
of $E\,.$ We prove that all points of this arc are boundary points of
the domain $D_k\,.$ Suppose that there exists a point
$z' \in \alpha \setminus \partial D_k\,.$ Because $z' \notin E,$
there exists a domain $D_{k'}$ such that $z' \in \partial D_{k'}, 1\le k' \le n\,.$
Choose sequences of points $z_m, z'_{m}, m=1,2,\dots$ and circular arcs
$\lambda_m\,, m=1,2,\dots$ in $U$ with $z_m, z'_{m}$  as end points and satisfying the
following condition. The points $z_m \in D_k, z_m \to z, m\to \infty\,;$
$ z'_{m} \in D_{k'}, z'_{m} \to z', m\to \infty\,,$ and the arcs
$\lambda_m\to T$ when $m \to \infty\,.$ This last convergence means that
upper bound of the distances of the points of the arc $\lambda_m$ to the
circle $T$ approaches $0$ when $m \to \infty\,.$ Because $D_k \cap D_{k'} =
\emptyset,$ we see that on every arc $\lambda_m$ there exists a point
$e_m \in E, m=1,2,3,\dots.$ The sequence of the points $\{e_m\}_{m=1}^{\infty}$
contains a subsequence converging to a point $e$ of the arc $\alpha\,.$
Because $E$ is closed, the point $e$ must be in $\alpha \cap E\,,$ a
contradiction with the choice of $\alpha\,.$ Let us denote by $\alpha_k$
the maximal open subarc of the circle $T\,,$ containing the arc $\alpha$
and not containing the points of the set $E\,.$ Clearly, the end points
of the arc $\alpha_k$ are contained in $E\,.$ The domain $D_k$ cannot
contain on its boundary two different arcs of the aforementioned type
(each corresponding to different points $z$) in view of the connectedness of
the set $E\,.$ We see that $  \alpha_k = (\partial D_k) \setminus E\,, k=1,2,\dots,n\,. $

We introduce the notation $B_k = D_k \cup \alpha_k \cup \{z: 1/ \overline{z} \in
D_k \}, k=1,\dots,n.$ From what was said above it follows that the sets
$B_k, k=1,\dots,n$ are simply connected pairwise nonintersecting domains on the
Riemann sphere $\overline{\mathbb C} \,.$ By Riemann's theorem, there exists
\begin{figure}[h]
 \includegraphics[scale=1]{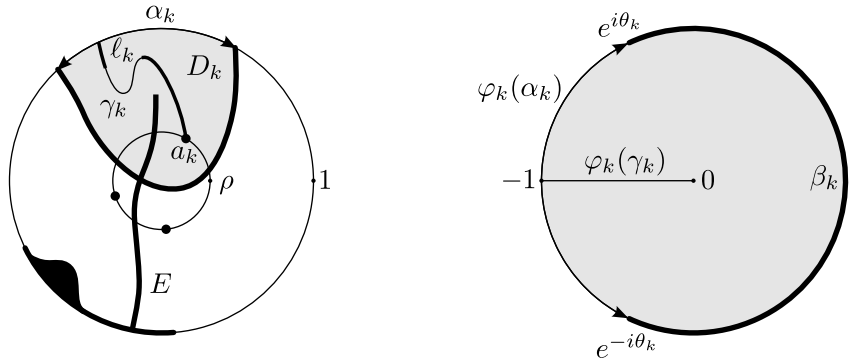}
\caption{Left: The continuum $E$ divides the unit disk $U$ into three domains.
Right: The image of the domain $D_k$ under $\varphi_k\,.$}
\end{figure}
a function $\varphi_k$ mapping conformally and univalently the corresponding
domain $B_k$ onto the domain $G_k$ which is the complement of the arc
$\beta_k:= \{  w: |w|=1, |{\rm arg} w| \le \theta_k \}, \varphi_k(a_k) = 0,$
$ \varphi_k'(a_k) > 0, \,\,\, k=1,\dots,n\,.$
On the basis of symmetry, the image of the domain $D_k$ is the disk $|w|<1\,,$
and further $\varphi_k(\alpha_k) = \{   w \in G_k: |w|=1\}$ and $\varphi_k(E \cap  \partial D_k) = \beta_k$. Because the harmonic measure is conformally invariant,
we have
$$
\omega_k = \omega(0, \beta_k, U_w)= \theta_k/ \pi, \qquad k=1,\dots,n\,.
$$
We denote by $\ell_k$ the union of a finite number of half-open arcs of the
analytic curve $\gamma_k:= \varphi_k^{-1}((-1,0])\,,$ having the following property:
for every $r, \rho\le r < 1\,,$ the circle $|z|= r$ intersects with $\ell_k$ only at one point. We set $z_k(r) = \{  z \in \ell_k: |z|= r\}\,.$ The construction of the arcs $\ell_k$
can be made in several ways. We choose one of them. It follows from the paper \cite{d1}
(cf. also \cite[p. 53]{d2}) that for every $r, \rho\le r <1\,,$ the following inequality
holds
\begin{equation} \label{d2ineq}
\sqrt[n]{ \Pi_{k=1}^n  r(B_k, z_k(r))} \le \frac{4 r}{n}\, ,
\end{equation}
and furthermore the equality in (\ref{d2ineq}) holds if and only if the domains
$B_k$ and the points $z_k(r)$ coincide with the domains $\{  z: |{\rm arg} z - 2 \pi(k-1)/n| < \pi/n \}$ and the points $r \, {\rm exp} (2 \pi (k-1)/n)$, respectively,
with a possible rotation around the origin. Here $r(B,z)$ stands for the inner radius
of the domain $B$ with respect to a point $z \in B\,, $ see \cite[Sec.1]{d2}.
The inequality (\ref{d2ineq}) was proved by the piecewise separating symmetrization method of
\cite{d2}. 

By the arithmetic-geometric mean inequality this implies that
$$
\frac{1}{n} \sum_{k=1}^n \frac{1}{ r(B_k, z_k(r))} \ge \frac{n}{ 4 r} \,.
$$
Integration yields the conclusion that
\begin{equation}  \label{my3}
\begin{aligned}
{} &  \displaystyle\frac{1}{n}  \sum_{k=1}^n  \int_{\gamma_k}  \displaystyle\frac{|dz|}{ r(B_k, z)}  \ge
\displaystyle\frac{1}{n}  \sum_{k=1}^n  \int_{\ell_k} \displaystyle\frac{|dz|}{ r(B_k, z)} \ge \\
{} & \displaystyle\frac{1}{n}  \sum_{k=1}^n  \int_{\rho}^1 \displaystyle\frac{d r}{ r(B_k, z_k(r))} \ge
\int_{\rho}^1 \displaystyle\frac{n}{4r}\, dr = - \displaystyle\frac{n}{4}  \log \rho  \,\, . \\
\end{aligned}
\end{equation}
On the other hand
\begin{equation} \label{my4}
\int_{\gamma_k} \frac{|dz|}{ r(B_k, z)} = \int_{-1}^0 \frac{dw}{ r(G_k, w)}
\end{equation}
for all $k=1,\dots,n\,.$ For the computation of the inner radius
$ r(G_k,w)$ we map the domain $G_k$ by a univalent conformal mapping onto the
unit disk $|\zeta|<1\,.$ It is easy to see that as such a map we can take a superposition
of the following functions:
$$
w_1= \frac{w-1}{w+1}\,, \qquad w_2 = -i w_1 {\rm ctg} \frac{\theta_k}{2}\,,\qquad \zeta= w_2- \sqrt{w_2^2-1}\,,
$$
where the last mapping is that branch of the function, inverse of the Joukowski
transformation, which maps $\infty$ to zero. Because under a conformal univalent
mapping the inner radius is multiplied by the modulus of the derivative at the
corresponding point, we have
$$
r(G_k,w) |w_1'(w)  {\rm ctg}  \frac{\theta_k}{2}  \zeta'(w_2)| = r(U_{\zeta}, \zeta) = 1- |\zeta|^2 \,.
$$
Simple calculations give from here
$$
r(G_k,w) = \frac{1-w}{ \sin ({\theta_k}/{2} )} \sqrt{w^2 - 2 w \cos {\theta_k} +1} \,.
$$
Inserting this value of the inner radius in the inequality (\ref{my4}), we obtain
$$
\int_{-1}^0 \frac{dw}{ r(G_k, w)} =  { \sin \frac{\theta_k}{2} } \int_{-1}^0 \frac{dw}{(1-w)\sqrt{w^2 - 2 w \cos {\theta_k} +1}} = \frac{1}{4}  \log \frac{1+ \sin ({\theta_k}/{2} )}{1- \sin ({\theta_k}/{2} )} \,.
$$
The computation of the integral is carried out either with standard integration methods or using tables such as \cite[p.94]{pbm}. Applying the relations (\ref{my3}) and (\ref{my4})
we arrive at the desired inequality (\ref{mainineq})\,.
\begin{figure}[h]
 \includegraphics[scale=1]{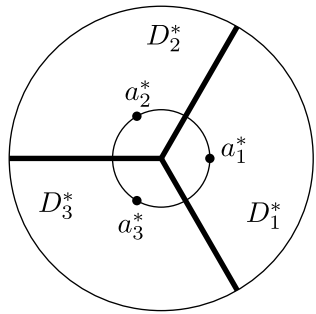}
\caption{Extremal configuration for $n=3\,.$}
\end{figure}

We now assume that in the inequality (\ref{mainineq}) we have equality. Then the equalities also hold in (\ref{my3}) and (\ref{d2ineq}) for all $r, \rho \le r < 1\, .$
From the equality in (\ref{my3}) it follows that the curves $\ell_k= \gamma_k$
are segments on radial rays, emanating from the origin, $k=1,\dots,n\,.$ Equality in (\ref{d2ineq}) shows that the angles between neighboring rays
are equal to $2\pi/n$ and the domains $B_k$
are sectors with angle $2\pi/n$ for which these rays are bisectors.
This observation implies that for some real
$\theta, a_k = a_k^*  e^{-i \theta}$ and $D_k = \{   z:   z  e^{i \theta} \in D_k^* \},
k=1,\dots,n\,.$ In this case the continuum $E$ necessarily contains the set
$E_{\theta}:=  \{ z: (z e^{i \theta})^n \in [-1,0]  \}\,.$ With a simple verification
we see that for $E= E_{\theta}$ in (\ref{mainineq}) we have the sign of equality.
It remains to observe that for other continua $E \supset E_{\theta}, E\neq E_{\theta}$
we have a strict inequality.
The theorem is proved.

\bigskip

{\bf Acknowledgements.} This work was completed during the visit of the first author
to the University of Turku, Finland. The authors are indebted to the referee, whose
comments were valuable.

\end{document}